\documentclass[reqno]{amsart}
\usepackage{amsthm,amsmath,amsfonts,amssymb,amscd,mathrsfs,graphics}
\usepackage{txfonts}
\usepackage{hyperref,supertabular}

\newcommand{\C}{{\mathbb{C}}}

\newcommand{\D}{{\mathscr D}}

\newcommand{\F}{{\mathscr F}}

\newcommand{\LL}{{\mathscr L}} 

\newcommand{\MM}{\overline{\M}} 
\newcommand{\M}{{\mathscr M}}

\newcommand{\Q}{{\mathbb{Q}}}

\newcommand{\Res}{{\operatorname{Res}}} 

\newcommand{\W}{{\mathscr W}}
\newcommand{\Z}{{\mathbb{Z}}}
\newcommand{\aut}{\operatorname{Aut}}
\newcommand{\balpha}{\boldsymbol{\alpha}}

\newcommand{\be}{\mathbf{e}}
\newcommand{\bgamma}{\boldsymbol{\gamma}} 

\newcommand{\bone}{\be}

\newcommand{\cC}{\mathscr C} 

\newcommand{\chat}{\hat{c}}
\newcommand{\ch}{{\mathscr H}}

\newcommand{\diag}{\operatorname{diag}}

\newcommand{\ga}{{\gamma}}

\newcommand{\grp}{G_W} 

\newcommand{\hGamma}{\widehat{\Gamma}}

\newcommand{\ind}{\operatorname{index}}

\newcommand{\Jac}{\operatorname{Jac}} 

\newcommand{\milnor}{\mathscr{Q}}

\newcommand{\restr}[2]{\left.#1\right|_{#2}} 

\newcommand{\st}{st} 

\newcommand{\unit}{\boldsymbol{1}}

\newcommand{\wit}{{\D}} 

\renewcommand{\hom}{\operatorname{Hom}}

\newtheorem{thm}{Theorem}[section]
\newtheorem{lm}[thm]{Lemma}
\newtheorem{prop}[thm]{Proposition}
\newtheorem{crl}[thm]{Corollary}
\newtheorem{conj}[thm]{Conjecture}

\theoremstyle{definition}
\newtheorem{rem}[thm]{Remark}

\newtheorem{df}[thm]{Definition}
\newtheorem{ex}[thm]{Example}
\theoremstyle{remark}
 
\title{Quantum ring of singularity $X^{p}+XY^{q}$}
\author{Huijun Fan}
\thanks{Partially Supported by NSFC 10401001, NSFC 10321001, and NSFC 10631050}
\address{School of Mathematical Sciences, Peking University, Beijing 100871, China}
\email{fanhj@math.pku.edu.cn}
\author{Yefeng Shen}
\address{School of Mathematical Sciences, Peking University, Beijing 100871, China}
\email{yfschen@gmail.com}
\begin{document}
\maketitle
\date{\today}

\begin{abstract}
In this paper, we will prove that the quantum ring of the
quasi-homogeneous polynomial $X^{p}+XY^{q}(p\ge 2,q>1)$ with some
admissible symmetry group $G$ defined by Fan-Jarvis-Ruan-Witten
theory is isomorphic to the Milnor ring of its mirror dual
polynomial $X^{p}Y+Y^{q}$. We will construct an concrete isomorphism
between them. The construction is a little bit different in case
$(p-1,q)=1$ and case $(p-1,q)=d>1$. Some other problems including
the correspondence between the pairings of both Frobenius algebras
has also been discussed.
\end{abstract}

\section{Introduction}

Let $(X,x)$ be an isolated complete intersection singularity of
dimension $N-1$. This means that $X$ is isomorphic to the fibre
$(f^{-1}(0),0)$ of an analytic map-germ $f:(\C^{N+k-1},0)\to
(\C^k,0)$, and $x\in X$ is an isolated singular point of $X$. In
particular, if $k=1$, $(X,x)$ is called a hypersurface singularity.
The study of the singularity was initiated by H. Whitney, R. Thom
and later developed by V. Arnold, K. Saito and many other
mathematicians during 60-80 years (See
\cite{AGV},\cite{S},\cite{ST},\cite{He}). The classification problem
is the central topic in the singularity theory. Many geometrical and
topological invariants were introduced to describe the behavior of
the singularity, for instance, the Milnor ring, intersection matrix,
Gauss-Manin system, periodic map and etc. The singularity theory has
tight connection with many fields in mathematics, like differential
equations, function theory and symplectic geometry.

Recently, in the papers \cite{FJR1,FJR2,FJR3} the first author and
his cooperators has constructed a quantum theory for the
hypersurface singularity if the singularity is given by a
non-degenerate quasi-homogeneous polynomial $W$. The start point of
their work is Witten's work \cite{Wi2} on the $r$-spin curves, where
Witten wanted to generalize the Witten-Kontsevich theorem
\cite{Wi1,K} to the moduli problem of $r$-spin curves (See
\cite{AJ},\cite{J1,J2},\cite{JKV1,JKV2},\cite{PV} for the discussion
on the $r$-spin curves). Unlike in the $r$-spin case that the Witten
equation has only trivial solution, in the general $W$ case, for
example $D_n$ and $E_7$ cases, the Witten equation may have
nontrivial solutions which can't be ignored in the construction of
the virtual cycle $\left[\W_{g,k}\right]^{vir}$. The Witten equation
is defined on an orbifold curve and has the following form
$$\bar{\partial}u_i+\overline{\frac{\partial W}{\partial u_i}}=0,$$
where $u_i$ are sections of appropriate orbifold line bundles.

The Witten equation comes from the study of the Landau-Ginzburg (LG)
model in supersymmetric quantum field theory (See \cite{Mar}). It
can be viewed as a geometrical realization of the $N=2$
superconformal algebra. The other known model is the Nonlinear sigma
model which corresponds to the Gromov-Witten theory in symplectic
geometry. In simple case, the LG model is totally determined by a
superpotential, which is a quasi-homogeneous polynomial required by
supersymmetry. There are two possible ways to get the topological
field theories by twisting the LG model. They are called the LG A
model and LG B model.The LG B model has been studied extensively in
physics and mathematics. The mathematical theory of LG A model (See
\cite{GS1,GS2} for the physical explanation) is just the quantum
singularity theory constructed by Fan-Jarvis-Ruan. As pointed out in
\cite{IV}, the more appropriate model is orbifold LG model, which
should be "identical" to a Calabi-Yau sigma model by CY/LG
correspondence. Actually the state space of the quantum singularity
theory is a space of the dual forms of Lefschetz thimbles
orbifolding by the admissible symmetric group $G$ of the polynomial
$W$.

Once we determine the state space and obtain the virtual cycle
$\left[\W_{g,k}\right]^{vir}$, we can build up the quantum
invariants for the singularity. For instance, we can define the
correlators $\langle\tau_{l_1}(\alpha_{i_1}), \dots,
\tau_{l_n}(\alpha_{i_n})\rangle^{W,G}_g$ for $\alpha_{i_j}$ in the
state space $\ch_{W,G}$ and the cohomological field theory. All the
correlators can be assembled into a generating function
$${\D}_{W,G}=\exp(\sum_{g\geq 0} \hbar^{2g-2}{\F}_{g, W,G}),$$
where
$${\F}_{g, W, G}=\sum_{k\geq 0} \langle\tau_{l_1}(\alpha_{i_1}), \dots,
\tau_{l_n}(\alpha_{i_n})\rangle^{W,G}_g\frac{t^{l_1}_{i_1} \cdots
t^{l_n}_{i_s}}{n!}$$ is the genus-$g$ generating function.

So computing out those quantum invariants becomes an important issue
to understand the singularity. Because of the Mirror symmetry
phenomena between the dual singularities (See \cite{Cl} and
references there), the quantum ring in the A model of the
singularity $W$ should be isomorphic to the Milnor ring in the B
model of the dual singularity $\check{W}$ (See
\cite{IV},\cite{Ka1,Ka2,Ka3}). Furthermore, we have more strong
conjecture relating the generating function ${\D}_{W,G}$ and the
formal Givental's generating function. Let us say more about this
conjecture.

The genus $g$ Gromov-Witten potential function of one point is
$$
\F_g^{pt}:=\sum_{k\ge 0}\frac{1}{k!}\sum_{d_1,\cdots,d_k}\langle
\tau_{d_1}\cdots,\tau_{d_k}\rangle_g t_{d_1}\cdots t_{d_k},
$$
where
$$
\langle
\tau_{d_1}\cdots,\tau_{d_n}\rangle_g=\int_{\MM_{g,k}}\psi^{d_1}_1\cdots\psi_k^{d_k}.
$$
The Witten-Kontsevich generating function is $\D^{pt}=\exp(\sum_g
\hbar^{g-1}\F^{pt}_g)$.

Let $A$ be a finite index set having a distinguish element $1$.
Suppose that the $\Q$ vector space $Vect(A)$ generated by $A$ is
attached with a nondegenerate symmetric bilinear form $\eta$. The
formal genus $0$ GW potential is a power series $\F_0$ in variables
$t_{d,l},d\in \mathbb{N}, l\in A$,
$$
\F_0=\sum_{k\ge
0}\frac{1}{k!}\sum_{\stackrel{d_1,\cdots,d_k}{l_1\cdots,l_k}}\langle
\tau_{d_1,l_1}\cdots,\tau_{d_k,l_k}\rangle_0 t_{d_1,l_1}\cdots
t_{d_k,l_k},
$$
which satisfies the string equation (SE), the dilaton equation (DE)
and the topological recursion equation (TRR).

Let $\F_{pr}=\F|_{\{t_{d_k,l_k}|t_{d_k,l_k}=0,\text{for} d_k>0\}}$
be the primary potential, then $\F_{pr}$ satisfies the WDVV equation
and form a Frobenius manifold. $\F_0$ is called semi-simple of rank
$\mu$ if $|A|=\mu$ and the algebra structure on $Vect(A)$ is
semi-simple for generic $t_{0,l}$. In \cite{Gi}, Givental found that
there is a transitive action of the so-called \emph{twisted loop
group} on the set of all semi-simple genus 0 GW potential of rank
$\mu$. Hence given a semi-simple potential $\F_0$ of rank $\mu$
there is group element $R$ taking $k$ copies
${\F_0^{pt}\oplus\cdots\oplus \F_0^{pt}}$ to $\F_0$.

Using a method to quantize the quadratic functions (see \cite{Gi}),
Givental can quantize the group element $R$ to get an element
$\hat{R}(\hbar)$ in \emph{Givental's group}. $\hat{R}(\hbar)$ acts
on the $k$ copies of the tau-functions $\D^{pt}\oplus\cdots\oplus
\D^{pt}$ to get a power series $\D_{Giv}$ in $\hbar$. $\D_{Giv}$ can
be written in the form $\D_{Giv}=\exp(\sum_g \hbar^{g-1}\F_g)$. If
$\D_{Giv}$ is required to satisfies a homogeneity condition, then
$\D_{Giv}$ is uniquely defined and satisfy the SE, DE, TRR and the
Virasoro constraints.

If given a genus 0 GW potential of a projective manifold which is
semi-simple, Givental conjectured that the total GW potential is the
same to $\D_{Giv}$ constructed from the genus 0 GW potential. We
have the similar question in the quantum singularity theory. Let
$\D_{W,G}$ be the tau-function in our Landau-Ginzburg A model and
$\D_{0,W,G}$ be the genus 0 tau-function. If the Frobenius manifold
induced by $\D_{0,W,G}$ is semi-simple, then we can get the formal
tau-function $\D_{Giv,W,G}$.

\begin{conj}\label{conj1}
$\D_{W,G}=\D_{Giv,W,G}$
\end{conj}

This should be true by Teleman's theorem \cite{Te} if we can show
that $\D_{0,W,G}$ is semi-simple. To prove the semi-simple property,
it is natural to show the Frobenius manifold associated to the
singularity $W/G$ in the A model is isomorphic to the Saito's
Frobenius manifold of the dual singularity $\check{W}$ in the B
model, which is easy proved to be semi-simple. If the symmetry group
$G$ is chosen suitable, we should have the following problem
\begin{conj}\label{conj2}
$\D_{W,G}$ is identical to $\D_{Giv,\check{W}}$ under some Mirror
transformation.
\end{conj}

Conversely, since the construction of the quantum theory depends on
the choice of the admissible subgroup $G$ such that $\langle
J\rangle \le G\le G_W$ (See the definitions in Section 2), we can't
expect a mirror correspondence from the LG A model of the dual
singularity $\check{W}$ with the trivial symmetry group to the LG B
model of the singularity $W$. A further discussion will be appeared
in \cite{Kr}.

In \cite{FJR2}, the authors has calculated the quantum ring
structure of the ADE singularities. Moreover by computing the basic
$4$ point correlators, using the WDVV equation and the
reconstruction theorems, the authors has proved the above conjecture
and thus proved the generalized Witten conjecture for DE cases via
the conclusions in \cite{GM}.

ADE singularities are simple singularities according to Arnold's
classification and has very special properties. For instance $ADE$
singularities are self-Mirror which has been shown in \cite{FJR2}.

To prove the Conjecture \ref{conj2} for general singularity, one has
to compare the Frobenius manifolds in both sides. Even in the B side
it is difficult to calculate the primary potential of Saito's
Frobenius manifold associated to a singularity other than ADE
singularities. One can consider the singularities with modality no
less than $1$. The computation of the quantum ring of those
singularities has recently been done (See \cite{Kr},\cite{Pr}). On
the other hand, M. Noumi \cite{No} has considered the following type
singularities:
\begin{enumerate}
\item[(i)] $x_1^{p_1}+x_2^{p_2}+\cdots+x_N^{p_N}$

\item[(ii)] $x_1^{p_1}+x_1x_2^{p_2}+x_3^{p_3}+\cdots+x_N^{p_N}$.
\end{enumerate}
He has considered the Gauss-Manin system associated to the above
singularity. An important fact is that the flat coordinates on the
Frobenius manifold are the polynomials of the deformed coordinates
appeared in the miniversal deformation, and meanwhile the formula of
primary potential was given in \cite{NY}.

Since the Frobenius structure of the above singularities in either
side is the tensor product of the Frobenius strucures of the $A_r$
singularity and the singularity $x^p+xy^q$, it is natural for us
only to compute the primary potential functions of the singularity
$x^p+xy^q$ in A model and then compare it with Noumi-Yamada's
computation in B model. By WDVV equation, one may show that the
primary potential depends on the $2,3$ point correlators and some
basic $4$ correlators. We need only compare the $2,3$ point
correlators and some $4$ correlators in both sides. The computation
of the quantum invariants of $x^p$ and $x^p+xy^q$ is important,
since we can take the direct sum of those singularities to form a
Calabi-Yau singularity (whose central charge is positive integer).
Once we know the quantum invariants of the Calabi-Yau singularity,
then by CY/LG correspondence it is hopeful to get the Gromov-Witten
invariants of the Calabi-Yau hypersurface defined by the CY
singularity. Actually a computation has been done in \cite{CR} for
quintic three-fold.

In this paper, we will calculate the quantum ring structure of the
singularity $x^p+xy^q,p\ge 2, q>1$ and construct the explicit
isomorphism to the Milnor ring of the dual singularity $x^py+y^q$ in
Berglund-H\"ubsch sense (see \cite{BH}). In a subsequent paper, we
will do the difficult computation of the basic $4$ point correlators
and build the isomorphism between two Frobenius manifolds. This
paper is arranged as follows. Section 2 gives a simple description
of the Fan-Jarvis-Ruan theory and list some useful axioms. In
Section 3, we will discuss the singularity $x^p+xy^q$ in case that
$(p-1,q)=1$. In section 4, we will treat the case that
$(p-1,q)=d>1$. In addition, in section 2, we also write down the
equivalence of the pairing of the dual forms of Lefschetz thimbles
and the residue pairing in the Milnor ring. This is just the Mirror
symmetry between the $2$ point functions. Though this fact was
mentioned in \cite{FJR2} and appeared in physical literature (see
\cite{Ce}), it is seldom known by mathematicians and is deserved to
be written down.

The first author would like to thank Tyler Jarvis, Masatoshi Noumi,
Kentaro Hori, and Marc Krawitz for their helpful discussion on the
related problem. Both authors would like to thank Yongbin Ruan for
his helpful suggestion, comment and kind help for many years.

\section{The Fan-Jarvis-Ruan-Witten theory}

\subsection*{The classical singularity theory} A polynomial $W:
\C^N \to \C$ is called quasi-homogeneous if there are positive
integers $d, n_1, \dots, n_n$ such that $W(\lambda^{n_1} x_1, \dots,
\lambda^{n_N} x_n)=\lambda^d w(x_1, \dots, x_N).$ We define the
\emph{weight} (or \emph{charge}), of $x_i$ to be
$q_i:=\frac{n_i}{d}.$ We say $W$ is \emph{nondegenerate} if (1) the
choices of weights $q_i$ are unique, and (2) $W$ has a singularity
only at zero. There are many examples of non-degenerate
quasi-homogeneous singularities, including all the nondegenerate
homogeneous polynomials and the famous $ADE$-examples:

\begin{ex}
\
    \begin{description}
    \item[$A_n$] $W=x^{n+1}, \ n\geq 1;$ \glossary{An@$A_n$ & The singularity defined by $W=x^{n+1}$}
    \item[$D_n$] $W=x^{n-1}+xy^2, \ n\geq 4;$ \glossary{Dn@$D_n$ & The singularity defined by $W=x^{n-1}+xy^2$}
    \item[$E_6$] $W=x^3+y^4;$\glossary{E6@$E_6$ & The singularity defined by $W=x^{3}+y^4$}
    \item[$E_7$] $W=x^3+xy^3;$\glossary{E7@$E_7$ & The singularity defined by $W=x^{3}+xy^3$}
    \item[$E_8$] $W=x^3+y^5;$\glossary{E8@$E_8$ & The singularity defined by $W=x^{3}+y^5$}
    \end{description}
    \end{ex}
A classical invariant of the singularity is the \emph{local
algebra}, also known as the \emph{Milnor ring} or \emph{Chiral ring}
in physics:
\begin{equation} \milnor_W:=\C[x_1, \dots, x_N]/\Jac(W),
\end{equation}
 where $\Jac(W)$ is the
Jacobian ideal, generated by partial derivatives
$$\Jac(W):=\left(\frac{\partial W}{\partial
x_1}, \dots, \frac{\partial W}{\partial x_N}\right).$$

The degree of the monomial makes the local algebra become a graded
algebra. There is a unique highest-degree element
$\det\left(\frac{\partial^2 W}{\partial x_i\partial x_j}\right)$
with degree \begin{equation}\hat{c}_W=\sum_i (1-2q_i).
\end{equation} which is called the \emph{central charge} of $W$.

The dimension of the local algebra is called the \emph{Milnor
number} and is given by the formula
$$\mu=\prod_i\left(\frac{1}{q_i}-1\right).$$.

Let $S$ be a small ball centered at the origin of $\C^\mu$ and
consider the miniversal deformation $F(x,t)$ of $W$ such that
$F(x,0)=W$. We have the Milnor fibration $F: \C^N\times S \to
\C\times S$ given by $(x,t)\to (F(x,t),t)$. Assume that the critical
value of $F$ are in $\C_\delta\times S$, where $\C_\delta:=\{z\in
\C:||z||<\delta\}$. Let $z_0\in \partial \C_\delta$, then
$F^{-1}(z_0,t)\to t\in S$ is a fiber bundle. This induces the
homology bundle $H_{N-1}(F^{-1}(z_0,t),\Z)\to S$. For a generic $t$,
$F(x,t)$ is a holomorphic Morse function. A distinguished basis of
$H_{N-1}(F^{-1}(z_0,t),\Z)$ can be constructed from a system of
paths connecting $z_0$ to the critical values. A system of paths
$l_i: [0,1]\to \C_\delta$ connecting $z_0$ to critical values $z_i$
is called \emph{ distinguished} if
\begin{itemize}
\item[(1)] $l_i$ has no self intersection;
\item[(2)] $l_i, l_j$ has no intersection except $l_i(0)=l_j(0)=z_0$;
\item[(3)] the paths $l_1, \dots,
l_{\mu}$ are numbered in the same order in which they enter the
point $z_0$, counter-clockwise.
\end{itemize}

For each $l_i$, we can associate a homology class $\delta_i\in
H_{N-1}(F^{-1}(z_0,t),\Z)$ as a vanishing cycle along $l_i$. More
precisely, the neighborhood of the critical point of $z_i$ contains
a local vanishing cycle. Then $\delta_i$ is obtained by transporting
the local vanishing cycle to $z_0$ using the homotopy lifting
property. The cycle $\delta_i$ is unique up to the homotopy of $l_i$
as long as the homotopy does not pass another critical value. Now
$\delta_1, \dots, \delta_{\mu}$ defines a distinguished basis of
$H_{N-1}(F^{-1}(z_0,t),\Z)$. The different choice of the
distinguished system of paths gives different distinguished basis.
The transformation relation between two basis is described by the
Picard-Lefschetz transformation. The intersection matrix
$(\delta_i\circ \delta_j)$ is an invariant of the singuarity and is
used to classify the singularity. Except the vanishing cycles,
another closely related objects are \emph{Lefschetz thimbles}, which
are the generators of the relative homology classes
$H_N(\C^N,F^{-1}(z_0,t),\Z)$. The boundary homomorphism $\partial$
gives an isomorphism $\partial: H_N(\C^N,F^{-1}(z_0,t),\Z)\to
H_{N-1}(F^{-1}(z_0,t),\Z)$. Geometrically, a Lefschetz thimble
$\Delta_i$ is the union of the vanishing cycles along the path $l_i$
and we have $\partial \Delta_i=\delta_i$.

We can let the radius $\delta$ of $\C_\delta$ goes to $\infty$ and
take $z_0=-\infty$. The relative homology class becomes
$H_N(\C^N,(Re F)^{-1}((-\infty, -M),t),\Z)$ for large $M>0$. We
simply write $(Re F)^{-1}((-\infty, -M),t)$ as $F_t^{-\infty}$ and
$(Re F)^{-1}((M,+\infty),t)$ as $F_t^{+\infty}$. Now the Lefschetz
thimble $\Delta_i$ in $H_N(\C^N,F_t^{-\infty},\Z)$ is canonically
determined by the horizontal path from the critical value to
$-\infty$.

Unlike the intersection matrix of the vanishing cycles, there is a
non-degenerate intersection pairing
\begin{equation}I:H_N(\C^N, F_t^{-\infty}, \Z)\otimes H_N(\C^N,
     F_t^{+\infty}, \Z)\to \Z.
\end{equation}

This pairing is given by the intersection of the stable manifold and
the unstable manifold of the critical point and is preserved by the
parallel transportation via the Gauss-Manin connection. Naturally we
have the dual pairing (See \cite{FJR2}):
$$\eta: H^N(\C^N, F^{-\infty}_{t}, \C)\otimes H^N(\C^N,
      F^{\infty}_{t}, \C)\to \C$$.

\subsection*{The Quantum invariants of the singularity}

Let $G_W:=\aut(W)$ be the \emph{maximal diagonal symmetry group} of
$W$ consisting of the diagonal matrix $\gamma$ such that $W(\gamma
x)=W(x)$. $G_W$ always contains the the subgroup $\langle J
\rangle$, where $J=\diag(e^{2\pi iq_1},\cdots,e^{2\pi i q_N})$ is
the \emph{exponential grading element}. We can take any subgroup $G$
such that $\langle J\rangle\le G\le G_W$. Using the group $G$, we
can orbifold the space of Lefschetz thimbles. For any $\gamma\in G$,
let $\C_{\gamma}^N$\glossary{cngamma@$\C_{\gamma}^N$ & The set of
fixed points of $\gamma$ in $\C^N$} be the set of fixed points of
$\gamma$, let $N_{\gamma}$\glossary{Ngamma@$N_{\gamma}$ & The
complex dimension of the fixed-point locus $\C^N_{\gamma}$} denotes
its complex dimension, and let $W_\ga:=\restr{W}{\C_{\gamma}^N}$ be
the quasi-homogeneous singularity restricted to the fixed point
locus of $\ga$. According to the Lemma 3.2.1 in [FJR2], $0$ is the
only critical point of $W_\gamma$ and $G$ is the subgroup of
$\aut(W_\gamma)$.

\begin{df}
The \emph{$\gamma$-twisted sector} $\ch_{\gamma}$ of the state space
is defined as the $G$-invariant part of the middle-dimensional
relative cohomology for $W_{\gamma}$;  that is,
\begin{equation}\label{eq:defChGamma}
\ch_\gamma:=H^{N_\gamma}(\C_\gamma^{N_{\gamma}},W_\gamma^{\infty},\C)^{G}.
\end{equation}
\end{df}

\begin{df} Suppose that $\gamma=(e^{2\pi i \Theta^\gamma_1}, \dots, e^{2\pi i\Theta^\gamma_{N}})\in G$ for
rational numbers $0\leq \Theta^\gamma_i<1$. The \emph{degree
shifting number} is $\iota_{\gamma}:=\sum_i (\Theta^{\gamma}_i
 -q_i)$ and for a class $\alpha \in \ch_{\gamma}$, we have the
 definition of the degree
\begin{equation*}
\deg_\C(\alpha):=\deg_W(\alpha)/2:=\deg(\alpha)/2+\iota_{\gamma}.
\end{equation*}
\end{df}

The following proposition was proved in Proposition 3.2.4 in [FJR2].

\begin{prop}\label{prp:iotaRel}
For any $\gamma \in G_W$ we have the equalities
\begin{align}
&\iota_{\gamma}+\iota_{\gamma^{-1}} = \hat{c}_W -N_\gamma\nonumber\\
&\deg_\C(\alpha)+\deg_\C(\beta) = \hat{c}_W
\end{align}
for any $\alpha \in \ch_\gamma$ and $\beta\in \ch_{\gamma^{-1}}$.
\end{prop}

\begin{df} The \emph{state space} of the
singularity $W/G$ is defined as
\begin{equation*}
\ch_{W,G}=\bigoplus_{\gamma\in G}\ch_{\gamma}.\end{equation*}
\end{df}

The pairing in $\ch_{W,G}$ is defined as the direct sum of the
pairings $$\langle \, , \rangle
_{\gamma}:\ch_{\gamma}\otimes\ch_{\gamma^{-1}}\to \C$$, where
$\langle \, , \rangle _{\gamma}$ is just the pairing
$\eta(\cdot,\cdot)$ of the singularity $W_\gamma$.

The quantum invariants of the singularity $W/G$ are defined via the
construction of the virtual fundamental cycle $\left[\W_{g,k}(
\bgamma)\right]^{vir}$ (or $\left[\W(\Gamma)\right]^{vir}$). Let us
briefly describe the properties of these virtual fundamental cycle
and some axioms related to our computation in this paper. We only
consider the case $G=G_W$ or $\langle J\rangle$.

Given a non-degenerate quasi-homogeneous polynomial $W$, we can
define the $W$-structure on an orbicurve with genus $g$ and $k$
marked points. Roughly speaking, the $W$ structure on a orbicurve
$\cC$ is a choice of $N$ orbifold line bundles $\LL_1,\cdots, \LL_N$
satisfying some relations defined by the polynomial $W$. If a
$W$-structure exists on an orbicurve $\cC$, then there must have
\begin{equation}\label{eq:sel-rule}
\deg(|\LL_j|) = \left(q_j(2g - 2 + k)
-\sum^k_{l=1}\Theta_j^{\gamma_l} \right)\in \Z.
\end{equation}
Here $\gamma_l=(e^{2\pi i \Theta^{\gamma_l}_1}, \dots, e^{2\pi
i\Theta^{\gamma_l}_{N}})\in G_W$ gives the orbifold action of the
line bundles $\LL_i$ at the marked point $z_l$ and $|\LL_j|$ is the
resolved line bundles on the coarse curve of $\cC$. (see [FJR2] for
the detail definition of these structures).

The orbicurve with $W$-structure is called \emph{$W$-orbicurve}. The
stack of stable $W$-orbicurves forms the moduli space $\W_{g,k}$.
For any choice $\bgamma:=(\gamma_1, \dots, \gamma_k) \in \grp^k$ we
define $\W_{g,k}(\bgamma)\subseteq \W_{g,k}$ to be the open and
closed substack with orbifold decoration $\bgamma$. We call
$\bgamma$ the \emph{type} of any $W$-orbicurve in
$\W_{g,k}(\bgamma)$. $\W_{g,k}(\bgamma)$ is not empty iff the
condition (\ref{eq:sel-rule}) holds. Forgetting the $W$-structure
and the orbifold structure gives a morphism $$\st:\W_{g,k } \to
\MM_{g,k}.$$  The morphism $\st$ plays a role similar to that played
by the {stabilization} morphism of stable maps in symplectic
geometry. The following theorem is proved in Theorem 2.2.6 of
[FJR2].

\begin{thm}
For any nondegenerate, quasi-homogeneous polynomial $W$, the stack
$\W_{g,k}$ is a smooth, compact orbifold (Deligne-Mumford stack)
with projective coarse moduli.  In particular, the morphism
$st:\W_{g,k} \to \MM_{g,k}$ is flat, proper and quasi-finite (but
not representable).  \end{thm}

Moreover, one can consider the decorated dual graph $\Gamma$ of a
stable $W$-curve and obtain the moduli space $\W_{g,k}(\Gamma)$,
which is a closed substack of $\W_{g,k}(\bgamma)$.

Let $T(\Gamma)$ be the set of tails of the decorated graph $\Gamma$
and attach an element $\gamma_\tau\in G_W$ to each tail $\tau$. The
virtual cycle $\left[\W(\Gamma)\right]^{vir}$ was constructed in
these papers $[FJR2,FJR3]$ by studying the Witten equation and its
moduli problem. It was proved that the virtual cycle
$\left[\W(\Gamma)\right]^{vir}$ satisfies a series of axioms
analogous to the Kontsevich-Manin axiom system in symplectic
geometry. We only list those axioms that we mainly used in this
paper.

Set \begin{equation} D:=-\sum_{i=1}^N
\ind(\LL_i)=\chat_W(g-1)+\sum_{j=1}^k \iota_{\gamma_j}.
\end{equation}

\begin{thm}\label{thm:main}\

\begin{enumerate}
\item \textbf{Dimension:}\label{ax:dimension} The cycle
  $\left[\W(\Gamma)\right]^{vir}$ has degree
\begin{equation}\label{eq:dimension}
6g-6+2k-2D=2\left((\chat-3)(1-g) + k - \sum_{\tau\in T(\Gamma)}
\iota_{\tau}\right).
\end{equation}
So the cycle lies in $H_r(\W(\Gamma),\Q)\otimes \prod_{\tau \in
  T(\Gamma)}
H_{N_{\gamma_{\tau}}}(\C^N_{\gamma_{\tau}},W^{\infty}_{\gamma_\tau},
\Q),$ where
$$
r:=6g-6+2k -2D -\sum_{\tau\in T(\Gamma)}N_{\gamma_\tau} =
2\left((\hat{c}-3)(1-g)+k-\sum_{\tau\in
T(\Gamma)}\iota(\gamma_{\tau})- \sum_{\tau\in
T(\Gamma)}\frac{N_{\gamma_\tau}}{2}\right).
$$

\item \textbf{Degenerating connected graphs:} \label{ax:ConnGraphs} Let $\Gamma$ be a
connected, genus-$g$, stable, decorated $W$-graph. The cycles
$\left[\W(\Gamma)\right]^{vir}$ and
$\left[\W_{g,k}(\bgamma)\right]^{vir}$ are related by
\begin{equation}
 \left[\W(\Gamma)\right]^{vir}=
\tilde{i}^*\left[\W_{g,k}(\bgamma)\right]^{vir},
\end{equation}
where $\tilde{i} : \W(\Gamma) \to{}{} \W_{g,k}(\bgamma)$ is the
canonical inclusion map.

\item{\bf Concavity}:\label{ax:convex}

    Suppose that all tails of $\Gamma$ are Neveu-Schwarz.  If
$\pi_*\left(\bigoplus_{i=1}^t\LL_i\right)=0$, then the virtual cycle
is given by capping the top Chern class of the dual $\left(R^1 \pi_*
\left(\bigoplus_{i=1}^t\LL_i\right)\right)^*$ of the pushforward
with the usual fundamental cycle of the moduli space:
\begin{equation}\begin{split}
\left[\W(\Gamma)\right]^{vir}& = c_{top}\left(\left(R^1\pi_*\bigoplus_{i=1}^t\LL_i \right)^*\right) \cap \left[\W(\Gamma)\right]\\
& = (-1)^D c_{D}\left(R^1\pi_*\bigoplus_{i=1}^t\LL_i \right) \cap
\left[\W(\Gamma)\right].
\end{split}
\end{equation}
\item   {\bf  Index zero:} \label{ax:wittenmap} Suppose that  $\dim (\W(\Gamma))=0$
and all the decorations on tails  are Neveu-Schwarz.

If the pushforwards $\pi_* \left(\bigoplus\LL_i\right)$ and
$R^1\pi_* \left(\bigoplus \LL_i\right)$ are both vector bundles of
the same rank, then the virtual cycle is just the degree
$\deg(\wit)$ of the Witten map times the fundamental cycle:
$$\left[\W(\Gamma)\right]^{vir} = \deg(\wit)\left[\W(\Gamma)\right],$$

\item\textbf{Composition law:}\label{ax:cutting} Given any genus $g$
  decorated stable $W$-graph $\Gamma$ with $k$ tails, and given any
  edge $e$ of $\Gamma$, let $\hGamma$ denote the graph obtained by
  ``cutting'' the edge $e$ and replacing it with two unjoined tails
  $\tau_+$ and $\tau_-$ decorated with $\gamma_+$ and $\gamma_-$,
  respectively.

The fiber product
$$F:=\W(\hGamma)\times_{\W(\Gamma)} \W(\Gamma)$$
has morphisms
 $$ \W({\hGamma})\xleftarrow{q} F
\xrightarrow{pr_2}\W(\Gamma).$$

We have
\begin{equation}\label{eq:cutting}
\left\langle
\left[\W(\hGamma)\right]^{vir}\right\rangle_{\pm}=\frac{1}{\deg(q)}q_*pr_2^*\left(\left[\W(\Gamma)\right]^{vir}\right),
\end{equation}
where $\langle  \rangle_{\pm}$ is the map from
$$H_*(\W(\hGamma)\otimes\prod_{\tau \in T(\Gamma)} H_{N_{\gamma_{\tau}}}(\C^N_{\gamma_{\tau}},W^{\infty}_{\gamma_\tau}, \Q) \otimes  H_{N_{\gamma_{+}}}(\C^N_{\gamma_{+}},W^{\infty}_{\gamma_+}, \Q)\otimes H_{N_{\gamma_{-}}}(\C^N_{\gamma_{-}},W^{\infty}_{\gamma_-}, \Q)$$ to
$$ H_*(\W(\hGamma)\otimes\prod_{\tau \in T(\Gamma)} H_{N_{\gamma_{\tau}}}(\C^N_{\gamma_{\tau}},W^{\infty}_{\gamma_\tau}, \Q)$$ obtained by contracting  the last two factors via the pairing
$$\langle\, , \rangle: H_{N_{\gamma_{+}}}(\C^N_{\gamma_{+}},W^{\infty}_{\gamma_+}, \Q)\otimes H_{N_{\gamma_{-}}}(\C^N_{\gamma_{-}},W^{\infty}_{\gamma_-}, \Q) \to \Q$$.
\end{enumerate}
\end{thm}

\noindent\emph{Cohomological field theory.} For any homogeneous
elements $\balpha :=(\alpha_1, \dots, \alpha_k)$ with $\alpha_i \in
\ch_{\gamma_i}$, the map $\Lambda^W_{g,k} \in \hom(\ch_W^{\otimes
k}, H^*(\MM_{g,k}))$ is defined by
 \begin{equation}
\Lambda^W_{g,k}(\balpha) := \frac{|G|^g}{\deg(\st)}PD
\,\st_*\left(\left[\W_{g,k}(W,\bgamma)\right]^{vir} \cap
\prod_{i=1}^k \alpha_i \right),
\end{equation}
and then extend linearly to general elements of $\ch_W^{\otimes k}$.
Here, $PD$ is the Poincare duality map.

The following results were showed in [FJR2]:

\begin{thm}\label{thm:CohFT}

The collection $(\ch_W, \langle\, , \rangle^W, \{\Lambda^W_{g,k}\},
\bone_1)$ is a cohomological field theory with flat identity.

Moreover, if $W_1$ and $W_2$ are two singularities in distinct
variables, then the cohomological field theory arising from
$W_1+W_2$ is the tensor product of the cohomological field theories
arising from $W_1$ and $W_2$:
$$(\ch_{W_1+W_2},  \{\Lambda^{W_1+W_2}_{g,k}\}) = (\ch_{W_1}\otimes\ch_{W_2},
\{\Lambda^{W_1}_{g,k} \otimes \Lambda^{W_2}_{g,k}\}). $$
\end{thm}

\begin{crl}
The genus-zero theory defines a Frobenius manifold.
\end{crl}

The quantum invariants of the singularity $W/G_W$ consists of the
correlators defined as below:
\begin{df}\label{df:correlator}
Define correlators $$\langle \tau_{l_1}(\alpha_1), \dots,
\tau_{l_k}(\alpha_k)\rangle^W_g:=\int_{\left[\MM_{g,k}\right]}\Lambda^W_{g,k}(\alpha_1,
\dots, \alpha_k)\prod_{i=1}^k {\psi}^{l_i}_i,$$ where $\psi_i$ are
the canonical classes in the tautological ring of $\MM_{g,k}$.
\end{df}

For an admissible group $G$ such that $\langle J\rangle \le G\le
G_W$, we can also define the virtual cycle
$\left[\W_{g,k,G}\right]^{vir}$, the morphism $\Lambda^{W,G}_{g,k}$
and the correlators $\langle \tau_{l_1}(\alpha_1), \dots,
\tau_{l_k}(\alpha_k)\rangle^{W,G}_g$. See the discussion in [FJR2].

\noindent\emph{Quantum ring (Quantum cohomology group) of the
singularity.} The simplest quantum structure of a singularity is the
Frobenius algebra consisting of the state space, the metric and the
quantum multiplication $\bigstar$. The multiplication is given by
the genus $0$ $3$-point correlators:
\begin{equation}
\langle\alpha\bigstar\beta, \gamma \rangle=\langle
\tau_0(\alpha),\tau_0(\beta),\tau_0(\gamma)\rangle_0^{W,G}.
\end{equation}

To show the mirror symmetry between the LG A model of the
quasi-homogenous singularity $W$ and the LG B model of the dual
singularity $\check{W}$, the first step is to identify the
corresponding Frobenius algebra structures and the second step is to
compare the Frobenius manifold structure. When the Frobenius
manifold structures are identical, it is hopeful to construct the
mirror map between the A model theory: Fan-Jarvis-Ruan-Witten theory
and the B model theory: Saito-Givental's theory.

Let us write down the explicit correspondence of the metric in A
model and the metric in B model. The identification was mentioned in
[FJR2] but not explicitly written down.

In LG B model, the Frobenius algebra is the Milnor ring $\milnor_W$
with the residue pairing and the multiplication of the monomials.
For $f,g\in \milnor_W$, the residue pairing is non-degenerate and is
defined by

$$\langle  f,g\rangle=\Res_{x=0}\frac{fg\, dx_1\wedge\cdots\wedge dx_N}{ \frac{\partial
W}{\partial x_1} \cdots \frac{\partial W}{\partial x_N}}.
$$

Let $\phi_i, i=0,1,\cdots, \mu-1$ be the basis of $\milnor_W$. We
can also consider the miniversal deformation
$F(x,t)=F_t(x):=W+t_0\phi_0+\cdots+t_{\mu-1}\phi_{\mu-1}$ and the
deformed Milnor ring with residue pairing
$\langle\cdot,\cdot\rangle_t=Res_t$ at the point $t\in \C^\mu$.

In the A model side, the intersection pairing $I$ of the Lefschetz
thimbles has the dual pairing
 $$\eta: H^N(\C^N, F^{-\infty}_{t}, \C)\otimes H^N(\C^N,
      F^{\infty}_{t}, \C)\to \C$$.
The relative cohomology groups $H^N(\C^N, F^{\pm\infty}_{t}, \C)$
and the pairing can be described in forms and the integration of
forms on $\C^N$.

Let $\bar{\partial}_{F_t}:=\bar{\partial}+d F_t\wedge,
{\partial}_{F_t}:={\partial}+d \bar{F}_t\wedge$, and $A^{p,q}$ be
the set of $(p,q)$-forms on $\C^N$. Then one can show that the
spectral sequence of the double complex $(A^{*,*}, \bar{\partial},d
F_t\wedge)$ converges to the homology group of
$(A^{*,*},\bar{\partial}_{F_t})$, which is also isomorphic to the
Koszul complex $(\Omega^*, dF_t)$. We obtain the isomorphisms
\begin{equation}
H^N_{\bar{\partial}_{F_t}}\simeq \Omega^N/dF_t\wedge
d\Omega^{N-1}\simeq \LL_{F_t}.
\end{equation}
Let $\{1=:\phi_0(x),\cdots,\phi_{\mu-1}(x)\}$ be a $\C$-basis in the
Milnor ring $\LL_{F_t}$ and $\omega=dx^1\cdots dx^N$ be the
holomorphic volume form in $\C^N$. Then the above isomorphisms can
be given canonically:
$$
\phi_i(x)\omega+\bar{\partial}_{F_t}\eta_i\longleftrightarrow
\phi_i(x)\omega\longleftrightarrow \phi_i(x).
$$
One can also study the cohomology group $H^N_{{\partial}_{F_t}}$
which is isomorphic to $H^N_{\bar{\partial}_{F_t}}$ by a usual
conjugation.

It is known that we can choose a family of primitive $n$-forms
$\{\omega_i\}$ ($\{\bar{\omega}_i\}$) generating
$H^N_{{\partial}_{F_t}}$ ($H^N_{\bar{\partial}_{F_t}}$). Such forms
are called vacuum wave forms in physical literature (see \cite{Ce}).
Note that $\bar{\partial}_{F_t}\omega_i={\partial}_{F_t} \omega_i=0
$. The closed forms $\{e^{F_t+\bar{F}_t}\omega_i\}$
($\{e^{-(F_t+\bar{F}_t)}*\bar{\omega}_i\}$ ) form a basis of
$H^N(\C^N, F^{-\infty}_{t}, \C)$ ($H^N(\C^N,F^{\infty}_{t}, \C)$).
Let $\{\Delta_a^-\,a=1,\cdots,\mu\}$ be a basis of
$H_N(\C^N,F^{-\infty}_{t}, \C)$ and $\{\Delta_b^+,b=1,\cdots,\mu\}$
be a basis of $H_N(\C^N,F^{\infty}_{t}, \C)$. Define
$$
\Pi^i_a=(-1)^{N/2}(2\pi)^{-N/2}\int_{\Delta_a^-}e^{F_t+\bar{F}_t}\omega_i;\;\tilde{\Pi}^j_b=
(-1)^{N/2}(2\pi)^{-N/2}\int_{\Delta_b^+}e^{-(F_t+\bar{F}_t)}*\bar{\omega}_i.
$$

Now the poincare dual $PD(\Delta_b^+)$ lies in $H^N(\C^N,
F^{-\infty}_{t}, \C)$, and we assume that $PD(\Delta_b^+)=\sum_i c_i
e^{F_t+\bar{F}_t}\omega_i$. Define
$\eta_{ij}=(-1)^N(2\pi)^{-N}\int_{\C^n} \omega_i\wedge
*\bar{\omega}_j$. By the relation
$$
\int_{\C^n} PD(\Delta_b^+)\wedge
e^{-(F_t+\bar{F}_t)}*\bar{\omega}_j=\sum_i c_i \int_{\C^n}
\omega_i\wedge *\bar{\omega}_j=(-1)^N(2\pi)^{N}\sum_i c_i \eta_{ij},
$$
we have
$$
c_i=(-1)^N(2\pi)^{-N}\sum_j
\eta^{ij}\int_{\Delta^+_b}e^{-(F_t+\bar{F}_t)}*\bar{\omega}_j.
$$

So we can compute the intersection number
\begin{align*}
&I_{a^- b^+}=\#(\Delta^-_a\cap
\Delta^+_b)=\int_{\Delta_a^-}PD(\Delta_b^+)\\
&=(-1)^N(2\pi)^{-N}\sum_{i,j}\int_{\Delta_a^-}
 e^{F_t+\bar{F}_t}\omega_i \eta^{ij}\int_{\Delta^+_b}e^{-(F_t+\bar{F}_t)}*\bar{\omega}_j\\
&=\sum_{ij} \Pi^i_a \eta^{ij}\tilde{\Pi}^j_b.
\end{align*}

\begin{rem}The vacuum wave forms $\omega_i$ can be chosen in the form
\begin{equation}\label{eq:vacuum-form}
\omega_i=\phi_i dx^1\cdots dx^N+\bar{\partial}_{F_t}\eta_i
\end{equation}
\end{rem}

The following result can be found in [Ce].

\begin{prop}\label{prop:resi=integration} Let $\omega_i$ have the representation
(\ref{eq:vacuum-form}) and
$\eta_{ij}=(-1)^N(2\pi)^{-N}\int_{\C^N}\omega_i\wedge
*\bar{\omega}_j$, then
\begin{equation}
\eta_{ij}=J(\phi_i dx^1\wedge\cdots\wedge dx^N,\phi_j
dx^1\wedge\cdots\wedge dx^N),
\end{equation}
where
\begin{equation}
J(\phi_i dx^1\wedge\cdots\wedge dx^N,\phi_j dx^1\wedge\cdots\wedge
dx^N):=Res_{x=0}\left(\frac{\phi_i\phi_j dx^1\wedge\cdots\wedge
dx^N}{\frac{\partial F_t}{\partial x^1}\cdots\frac{\partial
F_t}{\partial x^N}}\right)
\end{equation}
is the pairing in $\Omega^N/dF_t\wedge d\Omega^{N-1}$.
\end{prop}

This proposition shows that to compute the pairing $\eta:H^N(\C^N,
F^{-\infty}_{t}, \C)\otimes H^N(\C^N,
   F^{\infty}_{t},
\C)\to \C$ we need only compute the residue pairing of the
corresponding polynomials. Since the residue pairing is well-defined
at $t=0$, we can naturally extend the pairing $\eta$ at $t\neq 0$ to
$t=0$ by identifying it with the residue pairing. The identification
will be preserved if we consider the $G$-invariant theory.

\section{Quantum ring of $X^{p}+XY^{q};(p-1,q)=1$}\label{sec:2}

\subsection{Basic calculation}

Consider the singularity $W=x^{p}+xy^{q}$ with the constraint
$(p-1,q)=1, p\geq 2, q>1$. In this case, the group $G = \langle J
\rangle \cong \mathbb{Z} / (pq)\mathbb{Z}$. Let $\xi = exp(
\frac{2\pi i}{pq})$, then J acts on $\LL_{W} \omega$ by $(\xi^{p},
\xi ^{p-1})$.

We have the computation:
$$
\begin{matrix}
q_{x}=\frac{1}{p}, & q_{y}  = \frac{p-1}{pq},& \chat_{W}  = \frac{2(p-1)(q-1)}{pq} \\
\Theta_{x}^{J} = \frac{1}{p}, & \Theta_{y}^{J} = \frac{p-1}{pq}.&
\end{matrix}
$$
It is easy to obtain the state space
$$
\ch_{W,G} = \langle y^{q-1}\be_{0}, \be_{k} |k\in \Lambda\rangle,
$$
where $\Lambda=\{i \mid 1\leq i \leq pq-1, p \nmid{i}\}$, $\be_{0} :
= dx\wedge dy \in H^{mid}(\C_{J^{0}}^{N},W_{J^{0}}^{\infty}, \Q)$,
and $\be_{k} : =\unit \in H^{mid}(\C_{J^{k}}^{N},
W_{J^{k}}^{\infty}, \Q)$.

The complex dimension of $\ch_{W,G}$ is $ pq+1-q$.

Denote by $\{r\}$ the fractional part of the real number $r$. We
have
$$
\Theta_{x}^{J^{k}} =
\{\frac{k}{p}\},\;\Theta_{y}^{J^{k}}=\{\frac{k(p-1)}{pq} \}
$$
and the transition number
$$
\iota_{J^{k}} = \Theta_{x}^{J^{k}} - q_{x}+ \Theta_{y}^{J^{k}} -
q_{y} = \{\frac{k}{p}\} + \{\frac{k(p-1)}{pq}\} + \frac{1-p-q}{pq}
$$.
For any $\alpha\in \ch_{J^k}$, using the degree formula
$\deg_\C(\alpha)=\deg_{W}(\alpha)/2 = deg(\alpha)/2 +
\iota_{\gamma}$ we obtain
\begin{align*}
&\deg_{\C} (y^{q-1}\be_{0}) =
1-\frac{(1-p-q)}{pq} =
\frac{(p-1)(q-1)}{pq} = \hat{c}_{W}/2\\
&\deg_\C\be_{k} = \{\frac{k}{p}\} + \{\frac{k(p-1)}{pq}\} +
\frac{1-p-q}{pq}.
\end{align*}

\subsection{Computation of the 3-correlators of genus 0}

For convenience, we will write $y^{q-1}\be_{0}$ as $\be_{0}$ if
there is no confusion, and define the set $\hat{\Lambda}:=\Lambda
\bigcup \{0\}$.

The computation of the genus zero, three point correlators $\langle
a\be_{i}, b\be_{j}, c\be_{k} \rangle_{0}^{W}$. can be divided into
four cases.

\subsubsection*{Case1: $i=j=k=0$} By dimension formula, we have
$$
\langle
y^{q-1}\be_0,y^{q-1}\be_{0},y^{q-1}\be_0\rangle_{0}^{W} =0
$$.

\subsubsection*{Case2: only one of {i,j,k} not equal to 0}
The only non-zero correlator is
$\langle\be_1,y^{q-1}\be_{0},y^{q-1}\be_{0}\rangle_0^W$. Its value
is the residue pairing of the element $y^{q-1}\be_{0}$ with itself,
which is $-\frac{1}{q}$. Hence, we have
$$
\eta_{0,0}=\langle\be_1,y^{q-1}\be_{0},y^{q-1}\be_{0}\rangle_0^W=-\frac{1}{q}
$$,
and $\eta^{0,0}=-q$.

\subsubsection*{Case 3: $ijk\neq0$.}
\begin{lm}
If $ijk\neq0$, then
$\langle\be_{i},\be_{j},\be_{k}\rangle_{0}^{W}\neq0$ if and only if
$i+j+k$ equals to $pq+1$ or $2pq+1$. Furthermore,
$\langle\be_{i},\be_{j},\be_{k}\rangle_{0}^{W}=1$ if and only if the
corresponding line bundles satisfy $deg|\LL_x| = deg|\LL_y| = -1$
and $\langle\be_{i},\be_{j},\be_{k}\rangle_{0}^{W}=-q$ if and only
if $deg|\LL_x| = -2$ and $deg|\LL_y|=0.$
\end{lm}

\begin{proof}
If $ijk\neq0$, $\Sigma_{t\in \{i,j,k\}}deg_{\C}(\be_{t}) =
\Sigma_{t\in \{i,j,k\}}( \{ \frac{t}{p}\} +  \{ \frac{t(p-1)}{pq}\})
+ \frac{3(1-p-q)}{pq}$.

In this case, the degrees of two orbifold line bundles are
\begin{align*}
\deg|\LL_x| &= \frac{1}{p} - \{\frac{i}{p}\} - \{\frac{j}{p}\} -
\{\frac{k}{p}\}\\
\deg|\LL_y| &= \frac{p-1}{pq} - \{\frac{i(p-1)}{pq}\} -
\{\frac{j(p-1)}{pq}\} - \{\frac{k(p-1)}{pq}\}
\end{align*}

By the dimension counting, $\langle a\be_{i}, b\be_{j}, c\be_{k}
\rangle_{0}^{W}$ will vanish unless $\Sigma_{t}deg_{\C}(\be_{t}) =
\hat{c}_{W}$. Hence there is
$$deg|\LL_x| + deg|\LL_y| = -2$$

Since the degree of the resolved line bundles are integers, this
shows that
$$(i+j+k)(p-1) \equiv p-1 \mod (pq).$$

Since $3 \leq i + j + k \leq 3pq - 3$ and $(p-1,pq)=1$, we must have
$i + j + k = pq+1, $or$  2pq+1 $. Therefore

$$deg|\LL_x| = \frac{1}{p} - \{\frac{i}{p}\} - \{\frac{j}{p}\}
- \{\frac{k}{p}\} <0.$$

Since $deg|\LL_x| + deg|\LL_y| = -2$, we have two possibilities:
\begin{itemize}
\item[(1)] $deg|\LL_x| = deg|\LL_y| = -1$. In this case, by the concave axiom, we
have $\langle \be_{i}, \be_{j}, \be_{k} \rangle_{0}^{W}=1$.
\item [(2)] $deg|\LL_x| = -2$ and $deg|\LL_y| = 0$. In the same way, we
have $\langle \be_{i}, \be_{j}, \be_{k} \rangle_{0}^{W} = -q$.
\end{itemize}
\end{proof}

\begin{crl}\label{cr:metric} The metric has the form
$$
\eta_{\alpha\beta}=\left\{
\begin{array}{ll}
1, &
  \text{if}\;\alpha + \beta = pq \\
-1/q,& \text{if} \alpha=\beta = 0,
\end{array}
\right.
$$
\end{crl}

\begin{proof} It is obvious since we have the relation $\eta_{\alpha\beta} = \langle \be_{1},
\be_{\alpha}, \be_{\beta} \rangle_{0}^{W}$.
\end{proof}

\begin{rem} We also have the following conclusions:
\begin{enumerate}

\item For fixed $i$ and $j$, there is at most one $k$ such that
$\langle \be_{i}, \be_{j}, \be_{k} \rangle_{0}^{W} \neq 0$.
\item If
$2\leq (i+j) \leq pq$, then there must have $k = pq+1-(i+j)$.
\item
if $(pq+2)\leq (i+j) \leq (2pq-2)$, then $k= 2pq+1-(i+j)$.

\end{enumerate}
\end{rem}

\subsubsection*{Case4: if only one of $i,j,k$ equals to $0$}

\begin{lm}
Suppose only $k=0$ in $\{i,j,k\}$, then
$\langle\be_i,\be_j,y^{q-1}\be_0\rangle_{0}^{W}\neq0$ if and only if
$i+j=pq+1$ and $\langle \be_i,\be_j,\be_i,\be_j\rangle^W_0\neq 0$,
and furthermore
$\langle\be_i,\be_j,y^{q-1}\be_0\rangle_{0}^{W}=\pm{1}$.
\end{lm}

\begin{proof}
Suppose only $k=0$, then
$\langle\be_i,\be_j,y^{q-1}\be_0\rangle_{0}^{W}$ will vanish unless
$\deg_\C(\be_i) + \deg_{\C}(\be_j) + \deg_{\C}(y^{q-1}\be_0) =
\hat{c}_W$. Since $\deg_{\C}(y^{q-1}\be_{0}) = \hat{c}_{W}/2$, then
this is equivalent to
$$ \deg_{\C}(\be_i) + \deg_{\C}(\be_j) = \frac{(p-1)(q-1)}{pq}$$

On the other hand, by the composition axiom, we have
\begin{equation}\label{eq:compo}
\langle \be_i,\be_j,\be_i,\be_j\rangle_{0}^{W} =
\Sigma_{\alpha,\beta}\langle \be_i, \be_j,
\be_{\alpha}\rangle_{0}^{W} \eta^{\alpha \beta} \langle \be_i,
\be_j, \be_{\beta}\rangle_{0}^{W} + (\langle \be_i,
\be_j,y^{q-1}\be_0\rangle_{0}^{W})^{2}\eta^{0,0}
\end{equation}

Denote $\LL_{x,i,j,i,j}$ and $\LL_{y,i,j,i,j}$ by the orbifold line
bundles corresponding to the G-decorated graph (i,j,i,j)
respectively, then
\begin{align*}
\deg|\LL_{x,i,j,i,j}|&=2q_{x}-2\Theta_{x}^{J^i}-2\Theta_{x}^{J^j}=\frac{2}{p}-2\{\frac{i}{p}\}-2\{\frac{j}{p}\}\\
\deg|\LL_{y,i,j,i,j}|&=2q_{y}-2\Theta_{y}^{J^i}-2\Theta_{y}^{J^j}=\frac{2p-2}{pq}-2\{\frac{i(p-1)}{pq}\}-2\{\frac{j(p-1)}{pq}\}.
\end{align*}

Now we consider three cases.

\begin{itemize}

\item[(1)]. Both $p$ and $q$ are odd.

In this case, at least one of $\eta_{\alpha\beta}$,$\langle \be_i,
\be_j, \be_{\alpha}\rangle_{0}^{W}$,$\langle \be_i, \be_j,
\be_{\beta}\rangle_{0}^{W}$ vanishes, thus
$\langle\be_i,\be_j,y^{q-1}\be_0\rangle_{0}^{W} \neq 0$ if and only
if $\langle \be_i,\be_j,\be_i,\be_j\rangle_{0}^{W} \neq 0$.

Since the number $deg|\LL_{x,i,j,i,j}|$ and $deg|\LL_{x,i,j,i,j}|$
are integers, the dimension formula gives that
$deg|\LL_{x,i,j,i,j}|+deg|\LL_{x,i,j,i,j}|=-2$.

Since
$$deg|\LL_{x,i,j,i,j}|=\frac{2}{p}-2\{\frac{i}{p}\}-2\{\frac{j}{p}\}<0$$
and
$$deg|\LL_{y,i,j,i,j}|=\frac{2p-2}{pq}-2\{\frac{i(p-1)}{pq}\}-2\{\frac{j(p-1)}{pq}\}<1$$
we must have either
$$(\frac{2}{p}-2\{\frac{i}{p}\}-2\{\frac{j}{p}\},\frac{2p-2}{pq}-2\{\frac{i(p-1)}{pq}\}-2\{\frac{j(p-1)}{pq}\})=(-1,-1)$$
or
$$(\frac{2}{p}-2\{\frac{i}{p}\}-2\{\frac{j}{p}\},\frac{2p-2}{pq}-2\{\frac{i(p-1)}{pq}\}-2\{\frac{j(p-1)}{pq}\})=(-2,0)$$

Because $p$ is odd, only the second case is possible, thus by the
index-zero axiom, we have
$\langle\be_i,\be_j,\be_i,\be_j\rangle_{0}^{W}=-q$. Moreover, the
equality
$\frac{2p-2}{pq}-2\{\frac{i(p-1)}{pq}\}-2\{\frac{j(p-1)}{pq}\}=0$
implies $i+j=pq+1$.

The inverse conclusion is easy to see.

Now it is easy to check that
$\langle\be_i,\be_{pq+1-i},y^{q-1}\be_0\rangle_{0}^{W}=\pm{1}$.

\item[(2)]. $p$ is even and $q$ is odd.

In this case, the first term on the right hand of (\ref{eq:compo})
is not zero if and only if $\alpha=\beta=\frac{pq}{2}$ and
$deg|\LL_{x,i,j,\alpha}|=deg|\LL_{y,i,j,\alpha}|=-1$. This implies
that $deg|\LL_{x,i,j,i,j}|=deg|\LL_{y,i,j,i,j}|=-1$, i.e
$\langle\be_i,\be_j,\be_i,\be_j\rangle_{0}^{W}=1$. Thus
$\langle\be_i,\be_{j},y^{q-1}\be_0\rangle_{0}^{W}=0$. Here $i+j\neq
pq+1$.

If the first term on the right hand vanishes, then we have the same
conclusion by the same argument in (1).

\item[(3)]. $q$ is even.

Then by the assumption $(p-1,q)=1$, $p$ is even too. Hence
$\frac{pq}{2}\not\in \Lambda$, which concludes that
$\langle\be_i,\be_j,\be_i,\be_j\rangle_{0}^{W}\neq0$ if and only if
$\langle\be_i,\be_j,y^{q-1}\be_0\rangle_{0}^{W}\neq0$. So we have
$2(p-1)(1-i-j)\equiv 0 \mod (pq)$ since $deg|\LL_{y,i,j,i,j}|$ is an
integer. Then $i+j=\frac{pq}{2}+1,pq+1$ or $\frac{3pq}{2}+1$.

Assuming that $i+j=\frac{pq}{2}+1$ or $\frac{3pq}{2}+1$, we will
have
$$i+j\equiv 1 \mod (p)$$
 and
$$\frac{(i+j)(p-1)}{pq}\equiv \frac{p-1}{2}+\frac{p-1}{pq} \equiv
\frac{1}{2}+\frac{p-1}{pq} \mod (1).$$

Then $\{\frac{i}{p}\}+\{\frac{j}{p}\}=\frac{p+1}{p}$ and
$\{\frac{i(p-1)}{pq}\}+\{\frac{j(p-1)}{pq}\}=\frac{p-1}{pq}+\frac{1}{2}$
or $\frac{p-1}{pq}+\frac{3}{2}$.

So we have $ \deg_{\C}(\be_i) + \deg_{\C}(\be_j)=\Sigma_{i,j}(
\{\frac{i}{p}\} + \{\frac{i(p-1)}{pq}\} +
\frac{1-p-q}{pq})=\frac{(p-1)(q-1)}{pq}+1/2$ or
$\frac{(p-1)(q-1)}{pq}+3/2$ which contradicts with the degree
formula.

So only $i+j=pq+1$ is possible, then we can proceed as before to
reach the conclusion.
\end{itemize}

\end{proof}

\subsection{Generators and isomorphism}

Now, we will prove there exist two generators of the quantum ring
defined by Fan-Jarvis-Ruan-Witten theory. We need some preparation
to prove this fact.

\begin{lm}\label{lm: spec-inte} There exists a unique pair of integers $k$ and
$m$ in $\hat{\Lambda}$ satisfying the two conditions:
\begin{enumerate}
\item[(1)] $\deg_{\C}\be_k = \frac{(q-1)}{pq}$, $\deg_{\C}\be_m =
\frac{1}{q}$;

\item[(2)]
$$(k-1)(p-1) \equiv -1 \mod(pq)$$ and $$(m-2)(p-1) \equiv 1 \mod(pq)$$
\end{enumerate}
\end{lm}

\begin{proof}: Assume first that $p>2$. Because $(p-1,q)=1$, it is easy to see that the
congruence equation
$$(k-1)(p-1) \equiv -1 \mod(pq)$$ has a unique solution $k$ such that
$1\leq k \leq pq-1 $ and $p \nmid k$, this $k$ will satisfy $p\mid
(k-2)$, so $$\deg_{\C}\be_k =\{\frac{k}{p}\} + \{\frac{k(p-1)}{pq}\}
+ \frac{1-p-q}{pq}= \frac{2}{p} + \frac{p-2}{pq} - \frac{1}{p} -
\frac{p-1}{pq}= \frac{(q-1)}{pq}.$$

In the same way, we can find $\be_m$. Actually, there exists unique
$M$ , $1\leq M \leq (pq-1)$, such that
$$(p-1)M \equiv 1 \mod(pq).$$
Now we have $k = pq+1-M $ and set $m = M+2$.

If $p=2$, it is easy to get $M=1,k=0,m=3$ and $k\equiv pq+1-M$.
\end{proof}

In the following part, we still set $M=m-2$, where $m$ is the
special integer in Lemma \ref{lm: spec-inte}.

\begin{lm}
There is a bijective map $f$ between the set $\Delta:= \{(s,t)\in
\mathbb{Z}\bigoplus \mathbb{Z}\mid 0\leq s \leq p-2, 0 \leq t \leq
q-1\}$ and the set $\Lambda = \{ i \in \mathbb{Z} \mid 1 \leq i \leq
pq-1, p\nmid i \}.$
\end{lm}

\begin{proof}
We define a map $f : \Delta  \longrightarrow \Lambda$ as follows.

If there exists $i \in \Lambda $ such that $i \equiv 1+s(k-1)+t(m-1)
\mod (pq)$, then define $f(s,t) = i$. Since
$1+s(k-1)+t(m-1)=1+s(pq-M)+t(M+1)=spq + (M+1)(t-s) + (s+1)$, $p \mid
spq + (M+1)(t-s)$, so $p \nmid 1+s(k-1)+t(m-1)$. This thows that $f$
is well defined on $\Delta$.

Moreover, if $f(s,t)=f(s',t')$, then $1+s(k-1)+t(m-1) \equiv
1+s'(k-1)+t'(m-1) \mod (pq)$. We have $s-s' \equiv (M+1)(s'+t'-s-t)
\mod (pq)$. This implies $p \mid (s-s')$ and we must have $s'=s$ and
$(t'-t)(M+1) \equiv 0 \mod (pq)$. Thus we have $t'=t$ and the map is
injective.

Finally, since the sets $\Delta$ and $\Lambda$ have the same
cardinality $(p-1)q$, The map $f$ is bijective.
\end{proof}

Now for each $i \in \Lambda$, we can identify $\be_i,
\be_{1+s(k-1)+t(m-1)}$ and $(s,t)$, where $f(s,t) = i$.

\begin{lm}
$(s,t)\bigstar(u,v)=(s+u,t+v)$ if $0 \leq s+u \leq p-2 $ and $0 \leq
t+v \leq q-1$.
\end{lm}

\begin{proof}

Let $i=f(s,t), j= f(u,v)$,then
 $$(s,t)\bigstar(u,v) = \be_{i} \bigstar
\be_{j} = \Sigma_{\alpha,\beta} \langle \be_{i}, \be_{j},
\be_{\alpha} \rangle_{0}^{W} \eta^{\alpha \beta}\be_{\beta}$$

Because $i+j \equiv 1+s(k-1)+t(m-1)+1+u(k-1)+v(m-1) \equiv
(M+1)(t+v-s-u)+(2+s+u) \mod (pq)$, then $i+j \neq 1 \mod(pq) $. Now
by discussions above, we know there exists at most one element
$\alpha \in \Lambda$ such that $\langle \be_{i}, \be_{j},
\be_{\alpha} \rangle_{0}^{W} \neq 0$, and $\alpha \equiv
(pq+1)-(1+s(k-1)+t(m-1)+1+u(k-1)+v(m-1)) \equiv
-(M+1)(t+v-s-u)-(1+s+u) \mod (pq)$.

So we have $deg|\LL_x| = \frac{1}{p} - \{\frac{i}{p}\} -
\{\frac{j}{p}\} - \{\frac{\alpha}{p}\} = \frac{1}{p} - \frac{1+s}{p}
- \frac{1+u}{p} - \frac{p-1-s-u}{p} = -1$. Thus $deg|\LL_y| = -1$
also and $\langle \be_{i}, \be_{j}, \be_{\alpha} \rangle_{0}^{W} =
1$.

Now $(s,t)\bigstar(u,v)= \be_{i} \bigstar \be_{j} = \be_{pq -
\alpha} = \be_{1+(s+u)(k-1)+(t+v)(m-1)} = (s+u,t+v).$

\end{proof}

\begin{lm}
$(p-2,0)\bigstar(1,0)=\mp{q}\cdot y^{q-1}\be_0$
\end{lm}

\begin{proof}
$(p-2,0)\bigstar(1,0) = \be_{1+(p-2)(k-1)}\bigstar\be_k =
\Sigma_{\alpha\beta}\langle \be_{1+(p-2)(k-1)}, \be_k, \be_{\alpha}
\rangle_{0}^{W}\eta^{\alpha \beta}\be_{\beta} $. Since
$\be_{1+(p-2)(k-1)}=\be_M $ and $M+k=pq+1$, we have
 $\langle \be_{1+(p-2)(k-1)}, \be_k,
\be_{\alpha} \rangle_{0}^{W}\neq 0$ if and only if $\alpha = 0$.

Now we have $(p-2,0)\bigstar(1,0) = \langle \be_{M}, \be_k,
y^{q-1}\be_{0} \rangle_{0}^{W}\eta^{0,0}y^{q-1}\be_{0}=\mp{q}\cdot
y^{q-1}\be_{0}$.
\end{proof}

Define $(p-1,0)=\mp{q}\cdot y^{q-1}\be_{0}$. Then we have the
representation $(p-2,0)\bigstar(1,0)=(p-1,0)$ and
$(p-1-s,0)\bigstar(s,0)=(p-1,0)$ for any $0\leq s \leq p-2$.

\begin{lm}
$(p-1,0)\bigstar(0,1) = 0$.
\end{lm}

\begin{proof}
$(p-1,0)\bigstar(0,1)=\mp{q}\cdot y^{q-1}\be_{0}\bigstar \be_m
=\mp{q}{\Sigma_{\alpha\beta}\langle y^{q-1}\be_{0}, \be_m,
\be_{\alpha} \rangle_{0}^{W}\eta^{\alpha \beta}\be_{\beta}}$.

Now $\langle y^{q-1}\be_{0}, \be_m, \be_{\alpha} \rangle_{0}^{W}
\neq 0$ only if $m+\alpha = pq+1$. This implies $\alpha = pq+1-m =pq
- (M+1)$ and $p \mid \alpha$ which contradicts with $\alpha \in
\Lambda$.

Thus $\langle y^{q-1}\be_{0}, \be_m, \be_{\alpha} \rangle_{0}^{W} =
0$ for each $\alpha \in \Lambda$ and $y^{q-1}\be_{0}\bigstar \be_m =
0$.
\end{proof}

Define $(p,0)=(1,0)\bigstar(p-1,0)$.

\begin{lm}
$(p,0)+q(0,q-1)=0$.
\end{lm}

\begin{proof}
$(p,0) =(1,0)\bigstar(p-1,0) =\be_{k}\bigstar (\mp{q}\cdot
y^{q-1}\be_{0})
=\mp{q}\Sigma_{\alpha\beta}\langle\be_{k},y^{q-1}\be_{0},\be_{\alpha}
\rangle_{0}^{W}\eta^{\alpha \beta}\be_{\beta}
=\mp{q}\langle\be_{k},y^{q-1}\be_{0},\be_M
\rangle_{0}^{W}\eta^{M,pq-M}\be_{pq-M} =-q\be_{pq-M} =-q(0,q-1)$.
\end{proof}

\begin{lm}
$(s,t)\bigstar(u,v)= 0$ if $t+v \geq q$.
\end{lm}

\begin{proof}
$$(s,t)\bigstar(u,v)$$
$$=(s,0)\bigstar[(0,t)\bigstar(0,v)]\bigstar(u,0)$$
$$=(s,0)\bigstar[(0,t)\bigstar(0,q-1-t)\bigstar(0,1)\bigstar(0,v+t-q)]\bigstar(u,0)$$
$$=(s,0)\bigstar[(0,q-1)\bigstar(0,1)]\bigstar(0,v+t-q)\bigstar(u,0)$$
$$=(s,0)\bigstar[(-\frac{1}{q}(p,0)\bigstar(0,1)]\bigstar(0,v+t-q)\bigstar(u,0)$$
$$=(s,0)\bigstar0\bigstar(0,v+t-q)\bigstar(u,0)$$
$$=0.$$

\end{proof}

\begin{lm}
$(s,t)\bigstar(u,v)= 0$ if $s+u \geq p-1 $ and $t+v \neq 0$.
\end{lm}

\begin{proof}
$$(s,t)\bigstar(u,v)$$
$$=(0,t)\bigstar[(s,0)\bigstar(p-1-s,0)]\bigstar(s+u+1-p,v)$$
$$=(0,t)\bigstar(p-1,0)\bigstar(s+u+1-p,v)$$
$$=0.$$

\end{proof}

\begin{lm}
$(s,0)\bigstar(u,0)=-q(s+u-p,q-1)$ if $0 \leq {s,u} \leq p-1, p\leq
s+u\leq 2p-2 $.
\end{lm}

\begin{proof}
$(s,0)\bigstar(u,0) =(s,0)\bigstar(p-s,0)\bigstar(s+u-p,0)
=(p,0)\bigstar(s+u-p,0) =-q(0,q-1)\bigstar(s+u-p,0) =-q(s+u-p,q-1)$
\end{proof}

Now by the above lemmas, we obtain the following theorem.
\begin{thm}\label{thm: isom-1}
The two generators $\be_k $ and $\be_m$ in Lemma \ref{lm: spec-inte}
generate the quantum ring of $X^{p}+XY^{q}$, where $(p-1,q)=1$. The
multiplication is given by
\begin{enumerate}
\item $\mp{q}\cdot y^{q-1}\be_{0}=\be_k^{p-1}$;
\item $\be_{i} = \be_k^{s} \bigstar \be_m^{t}$, if $i \in \Lambda$ such that $f(s,t)=i$ for $(s,t) \in \Delta $.
\end{enumerate}
Moreover, we have two relations $\be_k^{p-1}\bigstar\be_m =
(p-1,0)\bigstar(0,1)=0$ and $\be_k^p + q\be_m^{q-1} = (p,0)+q(0,q-1)
= 0.$
\end{thm}

This theorem demonstrate the phenomenon of mirror symmetry between
two dual singularities:
\begin{crl} If $(p-1,q)=1$, then $\ch_{W,G} \cong \LL_{\check{W}}$, where
$\check{W}=X^{p}Y+Y^{q}$ is the dual singularity.
\end{crl}

\begin{proof}
Define the $\C$-algebra epimorphism $F: \C[X,Y]\longrightarrow
\ch_{W,G}$ such that $F(X)=\be_k$ and $F(Y)=\be_m$. Then
$F(X^{p-1}Y) = \be_k^{p-1}\bigstar \be_m = 0$ and $F(X^p + qY^{q-1})
= \be_k^{p} + q\be_m^{q-1} = 0$. Thus $X^{p-1}Y, X^p + qY^{q-1} \in
Ker(F)$, we have a $\C$-algebra epimorphism $\overline{F}:
\mathbb{C}[X,Y]/(pX^{p-1}Y, X^p + qY^{q-1}) \longrightarrow
\ch_{W,G}$. $\mathbb{C}[X,Y]/(pX^{p-1}Y, X^p + qY^{q-1})$ is just
the Milnor ring of the singularity $\check{W}=X^{p}Y+Y^{q}$. We have
$dim_{\C}\C[X,Y]/(pX^{p-1}Y, X^p + qY^{q-1}) = dim_\C\ch_{W,G}=
pq-q+1$. Those facts shows that $\overline{F}$ is a $\C$-algebra
isomorphism.
\end{proof}

\section{Quantum ring of $X^{p}+XY^{q};(p-1,q)=d>1$}

\subsection{Basic calculation}

We have the same fractional degrees and the central charge:

$$
q_{x}= \frac{1}{p},\;q_{y} = \frac{p-1}{pq},\;\hat{c}_{W}=
\frac{2(p-1)(q-1)}{pq}
$$.

Let $\xi = exp( \frac{2\pi i}{pq})$, and $\lambda$ acts on
$\mathcal{D}_{W} \omega$ by $(\xi^{-q}, \xi)$. Then $\lambda$
generates the maximal admissible abelian group $G=\Z / (pq)\Z$.

Now $\Theta_{x}^{J}=\frac{p-1}{p}$, and $\Theta_{y}^{J} =
\frac{1}{pq}$.

The G-invariant state space of the polynomial W is also:

$\ch_{W,G}=\langle y^{q-1}\be_{0}, \be_{k} |k\in \Lambda\rangle$,
where $\be_{0}, \be_{k}, \Lambda$ are defined as before.

The dimension is $\dim_\C\ch_{W,G}=pq+1-q$. We have the degree
computation:

\begin{align*}
&\Theta_{x}^{J^{k}}=\{\frac{p-k}{p}\},\;\Theta_{y}^{J^{k}}=\{
\frac{k}{pq} \}\\
&\iota_{J^{k}}=\Theta_{x}^{J^{k}}-q_{x}+ \Theta_{y}^{J^{k}} - q_{y}
= \{\frac{p-k}{p}\} + \{\frac{k}{pq}\} - \frac{p+q-1}{pq}\\
&\deg_{\C} (y^{q-1}\be_{0}) = 1-\frac{(p+q-1)}{pq} =
\frac{(p-1)(q-1)}{pq} = \hat{c}_{W}/2\\
&\deg_\C\be_{k} = \{\frac{p-k}{p}\} + \{\frac{k}{pq}\} +
\frac{1-p-q}{pq}.
\end{align*}

\begin{rem} Since $\deg\be_{p-1}=0$ in this case, So $\be_{p-1}$ will become the unit in the
quantum ring.
\end{rem}
\subsection{Computation of the $3$ point correlators in genus $0$}

The computation of the genus zero three point correlators $\langle
a\be_{i}, b\be_{j}, c\be_{k} \rangle_{0}^{W}$ is also divided into
four cases:

\subsubsection*{Case 1: $i=j=k=0$}

By dimension formula, we have
$$\langle
y^{q-1}\be_{0},y^{q-1}\be_{0},y^{q-1}\be_{0}\rangle_{0}^{W}=0.$$

\subsubsection*{Case 2: $j=k=0, i \in\Lambda$}

The only non-zero correlator is
$\langle\be_{p-1},y^{q-1}\be_{0},y^{q-1}\be_{0}\rangle_{0}^{W}$ and
$\langle\be_{p-1},\be_{j},\be_{k}\rangle_{0}^{W}=-\frac{1}{q}$.

\subsubsection*{Case 3: $i,j,k \in\Lambda$}

\begin{lm}
If $i,j,k \in\Lambda$, then there are two cases such that $\langle
a\be_{i}, b\be_{j}, c\be_{k} \rangle_{0}^{W}\neq0$.
\begin{itemize}
\item[(1)]
$\langle\be_{i},\be_{j},\be_{k}\rangle_{0}^{W}=-q$ if and only if
$i+j+k=p-1$
\item[(2)]
$\langle\be_{i},\be_{j},\be_{k}\rangle_{0}^{W}=1$ if and only if
$i+j+k=pq+p-1$

\end{itemize}
\end{lm}

\begin{proof} We have the computation:
\begin{align*}
&\Sigma_{i,j,k}\deg_{\C}(\be_{i}) = \Sigma_{i,j,k}( \{
\frac{p-i}{p}\} +  \{ \frac{i}{pq}\}) + \frac{3(1-p-q)}{pq}\\
&\deg|\LL_x| = \frac{1}{p} - \{\frac{p-i}{p}\} - \{\frac{p-j}{p}\} -
\{\frac{p-k}{p}\}\\
&\deg|\LL_y| = \frac{p-1}{pq} - \{\frac{i}{pq}\} - \{\frac{j}{pq}\}
-
\{\frac{k}{pq}\}\\
\end{align*}

If $\langle a\be_{i}, b\be_{j}, c\be_{k} \rangle_{0}^{W}\neq0$, then
by dimension formula, there is
$$
i+j+k \equiv p-1 \mod (pq).
$$

Thus $i+j+k=p-1,pq+p-1$ or $2pq+p-1$. But the later case implies
that $\Sigma_i\{\frac{p-i}{p}\}=\frac{1}{p}$, which is impossible.

If $i+j+k=p-1$, then $(deg|\LL_x|, deg|\LL_y|)=(-2,0)$ and
$\langle\be_{i},\be_{j},\be_{k}\rangle_{0}^{W}=-q$.

If $i+j+k=pq+p-1$, then $(deg|\LL_x|, deg|\LL_y|)=(-1,-1)$ and
$\langle\be_{i},\be_{j},\be_{k}\rangle_{0}^{W}=1$.
\end{proof}

\begin{rem}
The metric $\eta_{\alpha\beta}$ has the same form as in Corollary
\ref{cr:metric}.
\end{rem}

\subsubsection*{Case 4: $k=0,i,j\in\Lambda$}

\begin{lm}
$\langle\be_i,\be_j,y^{q-1}\be_0\rangle_{0}^{W}\neq0$ if and only if
$i+j=p-1$.

Moreover
$\langle\be_i,\be_{p-1-i},y^{q-1}\be_0\rangle_{0}^{W}=\pm{1}$.
\end{lm}

\begin{proof}

By the cutting formula, we have

\begin{align*}
\langle \be_i,\be_j,\be_i,\be_j\rangle_{0}^{W}
 &=\Sigma_{\alpha,\beta\in\Lambda}\langle \be_i, \be_j,
 \be_{\alpha}\rangle_{0}^{W} \eta^{\alpha \beta} \langle
\be_i, \be_j, \be_{\beta}\rangle_{0}^{W} +
(\langle\be_i,\be_j,y^{q-1}\be_0\rangle_{0}^{W})^{2}\eta^{0,0}\\
&=\Sigma_{l\in\Lambda}\langle \be_i, \be_j, \be_{l}\rangle_{0}^{W}
\eta^{l,pq-l} \langle \be_i, \be_j, \be_{pq-l}\rangle_{0}^{W} +
(\langle\be_i,\be_j,y^{q-1}\be_0\rangle_{0}^{W})^{2}\eta^{0,0}
\end{align*}
Now
\begin{align*}
&deg|\LL_{x,i,j,i,j}|=\frac{2}{p}-2\{\frac{p-i}{p}\}-2\{\frac{p-j}{p}\}\\
&deg|\LL_{y,i,j,i,j}|=\frac{2(p-1)}{pq}-2\{\frac{i}{pq}\}-2\{\frac{j}{pq}\}
\end{align*}

and the degree formula $\deg_{\C}(\be_i) + \deg_{\C}(\be_j) +
\deg_{\C}(y^{q-1}\be_0) =\hat{c}_W$ shows that
$deg|\LL_{x,i,j,i,j}|+deg|\LL_{y,i,j,i,j}|=-2$.

Since $-3<deg|\LL_{x,i,j,i,j}|<0$. Then $\langle
\be_i,\be_j,\be_i,\be_j\rangle_{0}^{W}\neq0$ if and only if
$(deg|\LL_{x,i,j,i,j}|,deg|\LL_{y,i,j,i,j}|)=(-2,0)$ or $(-1,-1)$.

On the other hand, $deg|\LL_{y,i,j,i,j}|=0$ or $-1$ implies
$i+j=p-1$ or $p-1+\frac{pq}{2}$ (if $2\mid(pq)$).

Then we have three cases:

\begin{itemize}
\item[(i)]

If $q$ is even, $\frac{pq}{2}$ is not belong to $\Lambda$, then
$\langle\be_i,\be_j,\be_i,\be_j\rangle_{0}^{W}=(\langle\be_i,\be_j,y^{q-1}\be_0\rangle_{0}^{W})^{2}\eta^{0,0}$.
Thus $\langle\be_i,\be_j,y^{q-1}\be_0\rangle_{0}^{W}\neq0$ implies
$i+j=p-1$ or $\frac{pq}{2}+p-1$.

If $i+j=p-1$,then
$deg|\LL_{x,i,j,i,j}|=\frac{2}{p}-2\{\frac{p-i}{p}\}-2\{\frac{p-j}{p}\}=-2$
and $deg|\LL_{y,i,j,i,j}|=0$,
$\langle\be_i,\be_j,\be_i,\be_j\rangle_{0}^{W}=-q$. If
$i+j=\frac{pq}{2}+p-1$, then $i+j\equiv -1 \mod (p)$. So
$deg|\LL_{x,i,j,i,j}|=\frac{2}{p}-2\{\frac{p-i}{p}\}-2\{\frac{p-j}{p}\}=-2$
and $deg|\LL_{y,i,j,i,j}|=-1$, which contradict with the fact that
$deg|\LL_{x,i,j,i,j}|+deg|\LL_{y,i,j,i,j}|=-2$;

Therefore $\langle\be_i,\be_j,y^{q-1}\be_0\rangle_{0}^{W}\neq0$ if
and only if $i+j=p-1$.
\item[(ii)] If $p$ is even and $q$ is odd, then $\frac{pq}{2}\in\Lambda$.

If $i+j=\frac{pq}{2}+p-1$,
$\langle\be_i,\be_j,\be_i,\be_j\rangle_{0}^{W}=1$, then we have
$$\Sigma_{l\in\Lambda}\langle \be_i, \be_j,
 \be_{l}\rangle_{0}^{W} \eta^{l,pq-l} \langle
\be_i, \be_j, \be_{pq-l}\rangle_{0}^{W} =\langle \be_i, \be_j,
 \be_{\frac{pq}{2}}\rangle_{0}^{W} \eta^{\frac{pq}{2},\frac{pq}{2}} \langle
\be_i, \be_j, \be_{\frac{pq}{2}}\rangle_{0}^{W} =1$$.

Thus $\langle\be_i,\be_j,y^{q-1}\be_0\rangle_{0}^{W}=0$.

If $i+j=p-1$, we have
$$
(\langle\be_i,\be_j,y^{q-1}\be_0\rangle_{0}^{W})^{2}\eta^{0,0}=\langle\be_i,\be_j,\be_i,\be_j\rangle_{0}^{W}
=-q.
$$
Thus $\langle\be_i,\be_j,y^{q-1}\be_0\rangle_{0}^{W}\neq0$ if
and only if $i+j=p-1$.

\item[(iii)] If both $p,q$ are odd, then $\frac{pq}{2}$ is not an integer.
Hence
$(\langle\be_i,\be_j,y^{q-1}\be_0\rangle_{0}^{W})^{2}\eta^{0,0}=\langle\be_i,\be_j,\be_i,\be_j\rangle_{0}^{W}
$. So $\langle\be_i,\be_j,y^{q-1}\be_0\rangle_{0}^{W}\neq0$ if and
only if $i+j=p-1$.
\end{itemize}
\end{proof}

\subsection{Generators and isomorphism}

\begin{lm}
There is a bijective map $g$ between the set $\Delta = \{(s,t)\in
\mathbb{Z}\bigoplus \mathbb{Z}\mid 0\leq s \leq p-2, 0 \leq t \leq
q-1\}$ and the set $\Lambda = \{ i \in \mathbb{Z} \mid 1 \leq i \leq
pq-1, p\nmid i \}.$
\end{lm}

\begin{proof}
We define a map $g : \Delta  \longrightarrow \Lambda$ as follows:

If there exists $i \in \Lambda $ such that $i \equiv p-1-s+tp \mod
(pq)$, then $g(s,t) = i$. Since $p \nmid tp-1-s$, this map is well
defined on the whole set $\Delta$.

Moreover, It is easy to verify that $g(s,t)=g(s',t')$ if and only if
$s'=s$ and $t'=t$. So the map is injective.

Finally, the cardinalities of the set $\Delta$ and the set $\Lambda$
are both equal to $(p-1)q$, we conclude that the map $g$ is
bijective.
\end{proof}

Now for any $i \in \Lambda$, we identify $\be_i$ and $(s,t)$ as the
same element if $g(s,t) = i$. Then $(1,0)=\be_{p-2}$ and
$(0,1)=\be_{2p-1}$. If we replace $\be_{k}$ and $\be_{m}$ in Theorem
\ref{thm: isom-1} by $\be_{p-2}$ and $\be_{2p-1}$ respectively, then
it is straightforward as in Section \ref{sec:2} to prove the
following theorem.

\begin{thm}
The two elements $\be_{p-2} $ and $\be_{2p-1}$ generate the quantum
ring of quasi-homogeneous polynomial $X^{p}+XY^{q}$ for $
(p-1,q)=d>1$. The multiplication is determined by the following
relations.
\begin{itemize}
\item $\be_{p-1}$ is the unit of this ring;
\item $\mp{q}\cdot y^{q-1}\be_{0}=\be_{p-2}^{p-1}$;
\item $\be_{i} = \be_{p-2}^{s} \bigstar \be_{2p-1}^{t}$ for each $i \in
\Lambda$, where $(s,t) \in \Delta $ such that $g(s,t)=i$,
\end{itemize}
and
\begin{itemize}
\item $\be_{p-2}^{p-1}\bigstar\be_{2p-1}=(p-1,0)\bigstar(0,1)=0;$
\item $\be_{p-2}^p + q\be_{2p-1}^{q-1} = (p,0)+q(0,q-1) =0.$
\end{itemize}
\end{thm}

\begin{crl} If $(p-1,q)=d>1$, and $G=\langle
\lambda\rangle$, then $\ch_{W,G}\cong \LL_{\check{W}}$. where
$\check{W}=X^{p}Y+Y^{q}$ is the dual singularity.
\end{crl}

\bibliographystyle{amsplain}

\begin{thebibliography}{JKV2}

\bibitem[AGV]{AGV} V. Arnold, A. Gusein-Zade and A. Varchenko,
\emph{Singularities of differential maps}, vol I, II, Monographs in
Mathematics, Birkhauser, Boston, 1985.

\bibitem[AJ]{AJ} D. Abramovich and T. Jarvis, \emph{Moduli of
twisted spin curves}.   Proc. of the Amer. Math. Soc.  \textbf{131}
(2003) no. 3, 685--699.

\bibitem[BH]{BH} P.~Berglund and T.~H\"ubsch,
\emph{A Generalized Construction of Mirror Manifolds}, Nucl. Phys. B \textbf{393} (1993) 377--391.

\bibitem[Ce]{Ce} S. Cecotti, $N=2$ Landau-Ginzburg VS. Calabi-Yau
$\sigma$-models: Non-perturbative Aspects, Int. J. Mod. Phys. A6
1749 (1991)

\bibitem[Cl]{Cl} P.~Clarke, \emph{Duality for toric Landau-Ginzburg models},
Preprint arXiv:0803.0447v1 [math.AG]

\bibitem[CR]{CR} A. ~Chiodo and Y. Ruan. \emph{Landau-Ginzburg/Calabi-Yau correspondence of quintic three-fold via
symplectic transformation}, Preprint arXiv: math/0812.4660v1.

\bibitem[FJR1]{FJR1} H. Fan, T. Jarvis and Y. Ruan,  \emph{Geometry and analysis of spin
equations}, Comm. Pure Applied Math 61(2008), 715-788.

\bibitem[FJR2]{FJR2} \bysame, \emph{The Witten equation, mirror symmetry
and quantum singularity theory}, Preprint arXiv:math/0712.4025?v3.

\bibitem[FJR3]{FJR3} \bysame, \emph{The Witten Equation and Its Virtual Fundamental Cycle}, book,
Preprint arXiv:math/0712.4025.

\bibitem[Gi]{Gi} \bysame. \emph{Gromov-Witten invariants and quantization of quadratic
Hamiltonians}, Dedicated to the memory of I. G. Petrovskii on the
occasion of his 100th anniversary, Mosc. J. 1 (2001), no. 4,
551-568, 645

\bibitem[GM]{GM} A. Givental and T. Milanov,  \emph{Simple singularities and integrable hierarchies},
The breadth of symplectic and Poisson geometry, Progr. Math., 232,
Birkhauser Boston, Boston, MA, (2005) 173--201.

\bibitem[GS1]{GS1} J. Guffin and E. Sharpe, \emph{A-twisted
Landau-Ginzburg models}. Preprint arXiv:0801.3836v1[hep-th]

\bibitem[GS2]{GS2} \bysame. \emph{A-twisted herotic Landau-Ginzburg
models.} Preprint arXiv:0801.3955v1[hep-th]

\bibitem[He]{He} C. Hertling, \emph{Frobenius manifolds and moduli
spaces for singularities}, Cambridge Univ. Press, Cambridge, 2002.

\bibitem[IV]{IV} K. Intriligator and C. Vafa, \emph{Landau-Ginzburg orbifolds}, Nuclear Phys.
B \textbf{339} (1990), no. 1, 95--120.

\bibitem[Ja1]{J1} T. J. Jarvis, \emph{Geometry of the moduli of higher
spin curves}, Inter. J. Math. \textbf{11} (2000), 637--663.

\bibitem[Ja2]{J2} \bysame, \emph{Torsion-free sheaves and moduli of
generalized spin curves}, Compositio Math. \textbf{110} (1998),
291--333.

\bibitem[JKV1]{JKV1} T. Jarvis, T. Kimura and A. Vaintrob, \emph{Moduli spaces
of higher spin curves and integrable hierarchies}, Compositio Math.
\textbf{126} (2001), 157--212.

\bibitem[JKV2]{JKV2} \bysame, \emph{Gravitational descendants and the
moduli space of higher spin curves},  Advances in algebraic geometry
motivated by physics (Lowell, MA, 2000), Contemp. Math.,
\textbf{276}, Amer. Math. Soc., Providence, RI, 2001, 167--177.

\bibitem[Ka1]{Ka1} R. Kaufmann, \emph{Singularities with symmetries, orbifold Frobenius algebras and mirror
symmetry},  Cont. Math. 403, 67-116.

\bibitem[Ka2]{Ka2} \bysame,  \emph{Orbifolding Frobenius algebras}, Internat. J. Math. \textbf{14} (2003), no. 6, 573--617.

\bibitem[Ka3]{Ka3} \bysame,  \emph{Orbifold Frobenius algebras, cobordisms and monodromies},
Orbifolds in mathematics and physics (Madison, WI, 2001), Contemp. Math., 310, Amer. Math. Soc.,
Providence, RI, 2002, pp. 135--161.

\bibitem[Ko]{K} M. Konstevich, \emph{Intersection theory on the moduli
space of curves and the matrix Airy function}, Comm. Math. Phys.
\textbf{164} (1992), 1--23.

\bibitem[Kr]{Kr} M.~Krawitz, \emph{FJRW-rings and mirror symmetry of singularities}, in preparation.

\bibitem[Mar]{Mar} E. Martinec, \emph{Criticality, Catastrophes,
and Compactifications}, in Physics and Mathematics of strings, ed L.
Brink, D. Friedan, and A. M. Polyakov.

\bibitem[No]{No} M. Noumi, \emph{Expansion of the solutions of a
Gauss-Manin system at a point of infinity}, Tokyo J. Math. Vol. 7,
No. 1, 1984

\bibitem[NY]{NY}  M. Noumi and Y. Yamada,
 \emph{Notes on the flat structures associated with simple and simply elliptic
  singularities},
  in "Integrable Systems and Algebraic Geometry," eds. M.-H.Saito, Y.Shimizu,
  and K. Ueno, World Scientific, (1998), 372--383.

\bibitem[Pr]{Pr} N. Priddis, P.~Acosta, M. Krawitz, N.~Wilde, and H.~Rathnakamura,
\emph{FJRW-ring and the Mirror Symmetry  for Bimodal Singularities.} Preprint.

\bibitem[PV]{PV} A. Polishchuk and A. Vaintrob, \emph{Algebraic construction
of Witten's top Chern class}, Advances in algebraic geometry motivated by physics (Lowell, MA, 2000), Contemp. Math., 276, Amer. Math. Soc., Providence, RI, 2001, pp. 229--249.

\bibitem[S]{S} K. Saito, \emph{Primitive forms for a universal
unfolding of a function with an isolated critical point}, J. Fac.
Sci. Univ. Tokyo Sect. IA 28 (1982)

\bibitem[ST]{ST} K. Saito, and Atsushi Takahashi, \emph{From primitive
forms to Frobenius manifolds}, preprint, 2008

\bibitem[Te]{Te} C. Teleman, The structure of 2D semi-simple field
theories, preprint 2007, Preprint arXiv: 0712.0160

\bibitem[Wi1]{Wi1} E. Witten, \emph{Two-dimensional gravity and
intersection theory on moduli space}, Surveys in Diff. Geom.
\textbf{1} (1991), 243--310.

\bibitem[Wi2]{Wi2} \bysame, \emph{Algebraic geometry associated with
matrix models of two-dimensional gravity}, Topological models in
modern mathematics (Stony Brook, NY, 1991), Publish or Perish,
Houston, TX, 1993, 235--269.

\end{thebibliography}

\providecommand{\bysame}{\leavevmode\hbox
to3em{\hrulefill}\thinspace}

\end{document}